\documentclass{article}

\usepackage[english]{babel}

\usepackage{amscd,amsmath,amsthm,amssymb,amsfonts,latexsym}
\theoremstyle{plain}
\newtheorem{theorem}{Theorem}

\newtheorem{definition}{Definition}

\def\bbbz{{\mathchoice {\hbox{$\sf\textstyle Z\kern-0.4em Z$}}
{\hbox{$\sf\textstyle Z\kern-0.4em Z$}} {\hbox{$\sf\scriptstyle
Z\kern-0.3em Z$}} {\hbox{$\sf\scriptscriptstyle Z\kern-0.2em Z$}}}}
\usepackage[dvips]{epsfig}

\theoremstyle{definition}
\newtheorem{remark}{Remark}
\newtheorem*{ack}{Acknowledgement}

\newcommand{\F}{{\mathbb{F}}}
\newcommand{\supp}{{\rm supp}}

\newcommand{\refT}[1]{Theorem~\ref{#1}}
\newcommand{\refS}[1]{Section~\ref{#1}}
\renewcommand\P{\operatorname{\mathbb P{}}}
\newcommand\floor[1]{\lfloor#1\rfloor}
\newcommand\ceil[1]{\lceil#1\rceil}
\newcommand\set[1]{\ensuremath{\{#1\}}}

\newcommand\bigpar[1]{\bigl(#1\bigr)}
\newcommand\Bigpar[1]{\Bigl(#1\Bigr)}

\newcommand\lrpar[1]{\left(#1\right)}

\newcommand\ma{s_\alpha}

\newcounter{CC}
\newcounter{cc}
\newcommand{\cc}{\stepcounter{cc}\ccx} 
\newcommand{\ccx}{c_{\arabic{cc}}}     
\newcommand\xfrac[2]{#1/#2}

\newcommand\qqw{^{-1/2}}
\newcommand\qq{^{1/2}}
\newcommand\ga{\alpha}
\newcommand\gao{{\alpha_0}}
\newcommand\mao{s_{\alpha_0}}
\newcommand\nao{n_{\alpha_0}}
\newcommand\gd{\delta}
\newcommand\na{n_\ga}
\newcommand\ta{t_\ga}
\newcommand\tao{t_{\alpha_0}}
\newcommand\ca{\mathcal A}
\newcommand\cax{\ca'}
\newcommand\ce{\mathcal E}
\newcommand\ea{\mathcal E_\ga}
\newcommand\eao{\mathcal E_\gao}
\newcommand\Bi{\operatorname{Bi}}
\newcommand\sumax{\sum_{\cax}}
\newcommand\EA{\bigwedge_{\ga\in\ca}\ea}
\newcommand\EAX{\bigwedge_{\ga\in\cax}\ea}
\newcommand\px{p_*}
\newcommand\eps{\epsilon}

\begin{document}

\title{On the Size of Identifying Codes in Binary Hypercubes}

\author{Svante Janson\\
Department of Mathematics\\ Uppsala University\\
P.O. Box 480\\SE-751 06 Uppsala, Sweden\\ and\\
 Tero Laihonen\thanks{Research
supported by the Academy of Finland under grant 111940.} \\
Department of Mathematics\\ University of Turku\\ FIN-20014 Turku,
Finland}
\date{18 April 2008}
\maketitle

\begin{abstract} In this paper, we consider
identifying codes in binary Hamming spaces $\F^n$, i.e., in binary
hypercubes. The concept of identifying codes was introduced by
Karpovsky, Chakrabarty and Levitin in 1998. Currently, the subject
forms a topic of its own with several possible applications, for
example, to sensor networks.

Let $C\subseteq \F^n$. For any $X\subseteq  \F^n$, denote by
$I_r(X)=I_r(C;X)$ the set of elements of $C$ within distance $r$ from at
least one $x\in X$. Now $C\subseteq \F^n$ is called an $(r,\leq
\ell)$-identifying code if the sets $I_r(X)$ are distinct for all
$X\subseteq \F^n$ of size at most $\ell$. Let us denote by
$M_r^{(\leq\ell)}(n)$ the smallest possible cardinality of an
$(r,\leq \ell)$-identifying code. In \cite{hlo3}, it is shown  for
$\ell=1$ that
$$\lim_{n\rightarrow\infty} \frac{1}{n}\log_2 M_r^{(\leq \ell)}(n)=1-h(\rho)$$
where $r=\lfloor \rho n \rfloor$, $\rho\in [0,1)$ and $h(x)$ is the
binary entropy function. In this paper, we prove that this result
holds for any fixed $\ell\geq 1$ when $\rho\in [0,1/2)$. We also
show that $M_r^{(\leq \ell)}(n)=O(n^{3/2})$ for every fixed $\ell$
and $r$ slightly less than $n/2$, and give an explicit construction
of small $(r,\leq 2)$-identifying codes for $r=\lfloor n/2 \rfloor
-1$.
\end{abstract}

\section{Introduction}

Let $\F=\{0,1\}$ be the binary field and denote by $\F^n$ the
$n$-fold Cartesian product of it, i.e. the Hamming space. We denote
by $A \ \triangle \ B$ the \emph{symmetric difference} $(A\setminus
B)\cup (B\setminus A)$ of two sets $A$ and $B$. The \emph{(Hamming)
distance} $d(x,y)$ between the vectors (called words) $x,y\in \F^n$
is the number of coordinate places in which they differ, i.e,
$x(i)\neq y(i)$ for $i=1,2,\dots, n$. The \emph{support} of
$x=(x(1),x(2),\dots,x(n)) \in \F^n$ is defined by $\supp(x)=\{i\mid
x(i)=1\}$. The \emph{complement} of a word $x\in \F^n$, denoted by
$\overline{x}$, is the word for which
$\supp(\overline{x})=\{1,2,\dots,n\}\setminus \supp(x)$. Denote by 0
the word where all the coordinates equal zero, and by 1 the all-one
word. Clearly $\overline{0}=1$. The \emph{(Hamming) weight} $w(x)$
of a word $x\in \F^n$ is defined by $w(x)=d(x,0)$. We say that $x$
$r$-\emph{covers}  $y$ if $d(x,y)\leq r$. The \emph{(Hamming) ball}
of radius $r$ centered at $x\in\F^n$ is
\[B_r(x)=\{y\in\F^n\mid d(x,y)\le r\}.\]
and its cardinality is denoted by $V(n,r)$.
 For  $X\subseteq \F^n$,
denote
\[
B_r(X)
=\bigcup_{x\in X} B_r(x)
=\{y\in\F^n\mid d(y,X)\le r\}
.
\]
We also use the notation
\[S_r(x)=\{y\in \F^n\mid d(x,y)=r\}.\]
A nonempty subset $C\subseteq \F^n$ is called a \emph{code} and its
elements are  \emph{codewords}. Let $C$ be a code  and $X\subseteq
\F^n$. We denote (the codeword $r$-neighbourhood  of $X$ by)
\[I_r(X)=I_r(C;X)=B_r(X)\cap C.\]
We write
for short $I_r(C;\{x_1,\ldots, x_k\})=I_r(x_1,\ldots, x_k)$.

\begin{definition} Let $r$ and $\ell$ be non-negative integers. A code
$C\subseteq \F^n$ is said to be $(r,\le \ell)$-\emph{identifying}
if for all $X,Y\subseteq \F^n$ such that $|X|\le \ell$, $|Y|\le
\ell$ and $X\ne Y$ we have
$$I_r(C;X)\ne I_r(C;Y).$$
\end{definition}

\medskip

The idea of the identifying codes is that given the set $I_r(X)$ we
can uniquely determine the set $X\subseteq \F^n$ as long as $|X|\leq
\ell$.

\smallskip

The seminal paper \cite{kcl} by Karpovsky, Chakrabarty and Levitin
initiated research in identifying codes, and it is nowadays a topic
of its own with different types of problems studied, see, e.g.,
\cite{BHL:ExtMinDenSqr},\cite{CGHLMM:LinAlgMin1-idTree},\cite{CHHLtwin-free},\cite{CHLext},\cite{fmmrs},\cite{GMS:IdCyc},\cite{M},\cite{RobRob};
for an updated bibliography of identifying codes see \cite{lowww}.
Originally, identifying codes were designed for finding
malfunctioning processors in multiprocessor systems (such as binary
hypercubes, i.e., binary Hamming spaces); in this application we
want to determine the set of malfunctioning processors $X$  of size
at most $\ell$ when the only information available is the set
$I_r(C; X)$ provided by the code $C$. A natural goal there is to use
identifying codes which is as small as possible. The theory of
identification can also be applied to sensor networks, see
\cite{rstu}.  Small identifying codes are needed for energy
conservation \cite{LT:disj}. For other applications we refer to
\cite{LTB-W:IdcodesSetVocer}.

The smallest possible cardinality of an $(r,\leq \ell)$-identifying
code in $\F^n$ is denoted by $M_r^{(\le\ell )}(n)$.

Let $h(x)=-x\log_2 x-(1-x)\log_2(1-x)$ be the binary entropy
function and $\rho\in [0,1)$ be a constant. Let further $r=\lfloor
\rho n \rfloor$. Honkala and Lobstein showed in \cite{hlo3} that,
when $\ell=1$, we have
\begin{equation}\label{aim}
\lim_{n\rightarrow\infty} \frac{1}{n}\log_2 M_r^{(\leq 1)}(n)=1-h(\rho).
\end{equation}
The lower bound that is part of \eqref{aim} comes from the simple
observation that if $C$ is an $(r,\le \ell)$-identifying code for
any $\ell\ge1$, then necessarily $B_r(C)=\F^n$ (otherwise there
would be a word $x\notin B_r(C)$ and then
$I_r(x)=\emptyset=I_r(\emptyset)$, so \set{x} and $\emptyset$ cannot
be distinguished by $C$) and also $|\F^n\setminus B_{n-r-1}(C)|\le1$
(otherwise there would be two words $x,y\notin B_{n-r-1}(C)$ and
then $I_r(\overline x)=C=I_r(\overline y)$, so \set{\overline x} and
\set{\overline y} cannot be distinguished by $C$); consequently, for
any $n, r, \ell\ge1$,
\begin{equation}\label{lb}
  \begin{split}
 M_r^{(\leq \ell)}(n)
\ge
 M_r^{(\leq 1)}(n)
&\ge
 \max\lrpar{
\frac{2^n}{|V(n,r)|},\,
\frac{2^n-1}{|V(n,n-r-1)|}}
\\
&=
 \max\lrpar{
\frac{2^n}{\sum_{i=0}^{r}\binom ni},\,
\frac{2^n-1}{\sum_{i=0}^{n-r-1}\binom ni}}
  \end{split}
\end{equation}
and the lower bound in \eqref{aim} follows from Stirling's formula.
Cf.\
\cite[Chapter 12]{chll}, \cite{bhl}, \cite{ELR:NewBounds},
\cite{hlo3} and \cite{kcl} for this and similar arguments and
related estimates.

Let us now consider any fixed $\ell>1$.
When $r=\lfloor \rho n \rfloor$,
we have by \eqref{aim} or
\eqref{lb}
the same lower bound as for $\ell=1$:
\begin{equation}\label{aiml}
\liminf_{n\rightarrow\infty} \frac{1}{n}\log_2 M_r^{(\leq \ell)}(n)
\ge 1-h(\rho).
\end{equation}
In the opposite direction,
it is shown in \cite{ELR:NewBounds} that
\begin{equation}\label{elr}
\limsup_{n\rightarrow\infty} \frac{1}{n}\log_2 M_r^{(\leq
\ell)}(n)\leq  1-(1-2\ell \rho)h\Bigpar{\frac{\rho}{1-2\ell \rho}}
\end{equation}
where $0\leq \rho\leq 1/(2\ell+1)$.
In this paper, we improve \eqref{elr} by showing that the lower bound
\eqref{aiml} is attained
for any fixed
$\ell\geq 1$ when $\rho\in [0,1/2)$.
(The proof is given in \refS{Spf1}.)

\begin{theorem}\label{T1}
  Let $\ell\ge1$ be fixed, let $\rho\in[0,1/2)$ and assume that
  $r/n\to\rho$. Then
\begin{equation*}
\lim_{n\rightarrow\infty} \frac{1}{n}\log_2 M_r^{(\leq
\ell)}(n)=1-h(\rho).
\end{equation*}
\end{theorem}

Furthermore, it is easy to see that when $\ell\ge2$,
unlike the case $\ell=1$, no
$(r,\leq \ell)$-identifying codes at all exist for
$r\geq \lfloor n/2\rfloor$. (This explains why we have to assume
$\rho<1/2$ in \refT{T1}.)

\begin{theorem}\label{Tnot}
If $n\geq 2$ and 
$r\geq \lfloor n/2\rfloor$,
then
there does not exist an $(r,\leq 2)$-identifying code
in $\F^n$.
\end{theorem}

\begin{proof}
If $r\geq \lfloor n/2\rfloor$, then
   $B_r(x)\cup B_r(\overline{x})=\F^n$
and thus  $I_r(x,\overline{x})=I_r(y,\overline{y})$ for any
$C\subseteq \F^n$ and $x,y\in \F^n$.
\end{proof}

We give this theorem mainly because the proof is so simple. In fact,
 it is proved in \cite{conf} that any $(r,\le \ell)$-identifying
code
 in $\F^n$ has to satisfy
\begin{equation}
  \label{sanna}
r\le \lfloor n/2\rfloor +2-\ell,
\end{equation}
which is slightly better than \refT{Tnot} when $\ell>3$.

Since $h(\rho)<1$ unless $\rho=1/2$, \refT{T1} implies that an
$(r,\le\ell)$-identifying code has to be exponentially large unless $r$
is close to $n/2$.
We give in \refS{Sconstr} an explicit construction
of a small
$(r,\leq 2)$-identifying code for
the largest possible $r$ permitted by \refT{Tnot}, viz.
$r=\lfloor n/2 \rfloor -1$.

\begin{theorem} \label{Tconstr}
Let $n\geq 2$.
There exists an $(r,\leq 2)$-identifying code in
$\F^n$ of size at most $n^3-n^2$
when $r=\lfloor n/2 \rfloor -1$.
\end{theorem}

For comparison,
it is shown in \cite{hlo3}
that for $\ell=1$ and $n\geq 3$,
$$M_{\lfloor \xfrac{n}{2}\rfloor}^{(\leq 1)}(n)\leq
\begin{cases}
\frac{n^2-n+2}{2}, & n \mbox{ odd, } \\
\frac{n^2-4}{2}, & n \mbox{ even. }
\end{cases}
$$

For $\ell>2$, we do not know any explicit constructions of small
$(r,\leq \ell)$-identifying codes in $\F^n$, but we can show the
existence of small such codes (even smaller than the one provided by
\refT{Tconstr}) for every $\ell\ge1$ when $r$ is a little smaller
than $n/2$. For $\ell=1$, there exist by the explicit estimate in
\cite[Corollary 13]{ELR:NewBounds} $(r,\leq 1)$-identifying codes in
$\F^n$ of size $O(n^{3/2})$ for every $r<n/2$ with
$r=n/2-O(\sqrt{n})$. Our next theorem, proved in \refS{Spf1}, yields
a bound of the same order
(although less explicit) for every fixed $\ell$ and certain $r$.

\begin{theorem}\label{T2}
Let $\ell\ge1$ be fixed and let $0<a<b$. Then there exist $n_0$ and
$A$ such that for every $n\ge n_0$ and $r$ with $n/2-b\sqrt n \le
r\le n/2-a\sqrt n$,
\begin{equation*}
M_r^{(\leq \ell)}(n) \le A n^{3/2}.
\end{equation*}
\end{theorem}

\begin{remark}
  This is not far from the best possible, since an $(r,\le
  \ell)$-identifying code $C$ in $\F^n$ trivially must satisfy
$\sum_{i=0}^\ell \binom{2^n}{i} \le 2^{|C|}$
and in particular
$2^n< 2^{|C|}$;
thus
$M_r^{(\leq \ell)}(n)>  n$ for  $r$ and $\ell\ge1$. (Moreover,
this argument yields $M_r^{(\leq \ell)}(n)\ge \ell n-O(1)$ for every
  fixed  $\ell\ge1$.)
\end{remark}

For $r$ closer to $n/2$, we can show a weaker result,
still with a polynomial bound. (This theorem too is proved in \refS{Spf1}.)
\begin{theorem}\label{T3}
Let $\ell\ge1$ be fixed and let $L$ be fixed with $L\ge 2^\ell$.
Then there exist $n_0$ and $A$ such that for every $n\ge n_0$ and $r$
with
$r=\lfloor n/2 \rfloor -L$,
\begin{equation*}
M_r^{(\leq \ell)}(n) \le A n^{2^{l-1}+1}.
\end{equation*}
\end{theorem}
For $\ell\ge3$, we do not know
the largest possible $r$ such that
there exists  an $(r,\leq \ell)$-identifying code
in $\F^n$,
but \refT{T3}
leaves only a small gap to the bounds
in \refT{Tnot} and \eqref{sanna}.

\section{Proofs of the main results}\label{Spf1}

Our non-constructive upper bounds in
Theorems \ref{T1}, \ref{T2} and \ref{T3}
are based on the following general theorem proven in
\cite{ELR:NewBounds}.
Let
$m_n(r,\ell)$ stand for the minimum of $|B_r(X)\ \triangle \
B_r(Y)|$ over any subsets $X,Y\subseteq \F^n$, $X\neq Y$ and $1\leq
|X|\leq \ell$ and $1\leq |Y|\leq \ell$.  Denote further by $N_\ell$
the number of (unordered) pairs $\{X, Y\}$ of subsets of $\F^n$ such
that $X\neq Y$ and $1\leq |X|\leq \ell$ and $1\leq |Y|\leq \ell$.

\begin{theorem}[\cite{ELR:NewBounds}]
\label{basic}
Let $r\geq 1$, $\ell\geq 1$ and $n\geq 1$.
Provided that $m_n(r,\ell)>0$, there exists an $(r,\leq \ell)$-identifying code
of size $K$ in $\F^n$ such that \[K\leq \left\lceil
\frac{2^n}{m_n(r,\ell)} \ln N_\ell\right\rceil +1.\]
\end{theorem}

Obviously,
\[N_\ell\leq \left(\sum_{i=1}^\ell \binom{2^n}{i} \right)^2
\le 2^{2n\ell}
\]
and thus \refT{basic} yields
\begin{equation}\label{bas}
M_r^{(\leq \ell)}(n) \le \frac{2^{n+1} \ell n}{m_n(r,\ell)} +2.
\end{equation}

It remains to estimate $m_n(r,\ell)$. Using probabilistic arguments,
we are able to show in Theorems \ref{T1m} and \ref{T2m} the
following crucial result: for fixed $\ell$ and $\eps>0$ or $a,b$,
there exists constants $n_0$ and $c>0$ such that for $n\geq n_0$ and
$r$ with $r\le(1/2-\eps)n$, or $n/2-b\sqrt n \le r\le n/2-a\sqrt
n$,
we have
\begin{equation} \label{estimate}
m_n(r,\ell)\geq c \binom{n}{r}.
\end{equation}
(By combing the methods of proofs below, it is possible to show that
this holds in the intermediate range of $r$ too, but we omit the details.)
For $\ell$ fixed and $r=\floor{n/2}-L$, with $L\ge 2^\ell$ fixed,
we prove in \refT{T3m} the slightly weaker
estimate
\begin{equation} \label{estimate2}
m_n(r,\ell)\geq c n^{-2^{\ell-1}}2^n.
\end{equation}
(We do not know whether \eqref{estimate} holds in this case too.)
Combining \eqref{bas},  \eqref{estimate}, \eqref{estimate2} and
standard estimates for binomial coefficients, see \cite[p.\
33]{chll}, we obtain Theorems \ref{T1}, \ref{T2} and \ref{T3}.

We prove the required estimates of $m_n(r,\ell)$ in the following
form. In applying the following results to obtain the bounds
\eqref{estimate} and \eqref{estimate2} on $m_n(r,\ell)$ just notice
that we can assume that there is $x\in X\setminus Y$ and $Y\subseteq
\{y_1,\dots,y_\ell \}$.

\begin{theorem} \label{T1m}
Let $\ell\geq 1$ be fixed.
For every $\varepsilon >0$ there is a constant $c>0$ and $n_0$ such
that for $n\geq n_0$ and any $\ell+1$ words $x$ and $y_1,\dots,
y_\ell$ in $\F^n$, with  $y_i\neq x$ for $i=1,\dots,\ell$, and every
$r$ with $0\leq r\leq (1/2-\varepsilon)n$, there exist at least $c
\binom{n}{r}$ words $z\in \F^n$ with $d(z,x)=r$ and $d(z,y_i)>r$ for
$i=1,\dots,\ell$.
\end{theorem}

\begin{theorem} \label{T2m}
Let $\ell\geq 1$ be fixed. For every $a,b>0$ there is a constant
$c>0$ and $n_0$ such that for $n\geq n_0$ and any $\ell+1$ words $x$
and $y_1,\dots, y_\ell$ in $\F^n$, with  $y_i\neq x$ for
$i=1,\dots,\ell$, and every $r$ with $n/2-b\sqrt n \le r\le
n/2-a\sqrt n$,
there exist at least $c n\qqw 2^n \ge
c\binom{n}{r}$ words $z\in \F^n$ with $d(z,x)=r$ and $d(z,y_i)>r$
for $i=1,\dots,\ell$.
\end{theorem}

\begin{theorem} \label{T3m}
Let $\ell\geq 1$ be fixed.
For every $L \ge2^\ell$ there is a constant $c>0$ and $n_0$ such
that for $n\geq n_0$ and any $\ell+1$ words $x$ and $y_1,\dots,
y_\ell$ in $\F^n$, with  $y_i\neq x$ for $i=1,\dots,\ell$, and
$r=\floor{n/2}-L$,
there exist at least $cn^{-2^{\ell-1}}2^n$ words $z\in \F^n$ with
$d(z,x)=r$ and $d(z,y_i)>r$ for
$i=1,\dots,\ell$.
\end{theorem}

The proofs of Theorems \ref{T1m}--\ref{T3m} are similar, although
some details differ. We begin with some common considerations.

By symmetry we may assume that $x=0$. Given $y_1,\dots, y_\ell$,
partition the index set $[n]=\{1,\dots,n\}$ into $2^\ell$ subsets
$A_\alpha$, indexed by $\alpha\in \F^\ell$, such that
$$A_\alpha=\{i\in [n]: y_j(i)=\alpha_j \mbox{ for }
j=1,\dots,\ell\}.$$

Let $z\in \F^n$ and let further $\ma=\ma(z)=|\{i\in
A_\alpha : z(i)=1\}|$. Then $d(z,x)=\sum_{\alpha}\ma$ and
$$d(z,y_j)=\sum_{\alpha:\alpha_j=0}\ma+\sum_{\alpha:\alpha_j=1}(|A_\alpha|-\ma)
=d(z,x)+\sum_{\alpha:\alpha_j=1}(|A_\alpha|-2\ma).$$
Hence, if $d(z,x)=r$, we need also
$$\sum_{\alpha:\alpha_j=1}(|A_\alpha|-2\ma)\geq 1$$
for each $j=1,\dots,\ell$; then $d(z,y_j)>d(z,x)=r$.

For simplicity, we consider only $z$ such that
$\ma<|A_\alpha|/2$ for every $\alpha$ such that $A_\alpha\neq
\emptyset$; we say that such $z$'s are \emph{good}.
Note that $\sum_{\alpha:\alpha_j=1}|A_\alpha|=d(x,y_j)\geq 1$ for
each $j$, so $A_\alpha\neq \emptyset$ for some $\alpha$ with
$\alpha_j=1$, and if $z$ is good, then
$\sum_{\alpha:\alpha_j=1}(|A_\alpha|-2\ma)>0$, and thus, as
shown above, we get $d(z,y_j)>d(z,x)$ for each $j$.
Thus, it suffices to show that the number of good words $z$ with
$d(z,x)=r$ is at least the given bounds in the theorems.

\begin{proof}[Proof of \refT{T1m}]
It now suffices to show that there exist $c$ and $n_0$ such that
for any choice of $n\geq n_0$, $x=0$, $y_1,\dots,y_\ell$ and $r$
with $0\leq r\leq (1/2-\varepsilon)n$, if $z$ is a random word with
$d(z,x)=r$, i.e., a random string of $r$ $1$'s and $n-r$ 0's, then
$$\P(z \mbox{ is good})\geq c.$$

Suppose that this is false for all $c$ and $n_0$. Then there exists
a sequence of such $(n,y_1,\dots,y_\ell, r)$, say
$n_\gamma,y_1^{(\gamma)},\dots,y_\ell^{(\gamma)}\in \F^{n_{\gamma}}$
and $r_\gamma$, $\gamma=1,2,\dots$, such that
$n_\gamma\to\infty$ and
if $z\in
\F^{n_{\gamma}}$ is a random string with $r_\gamma$ $1$'s, then
$$\P(z \mbox{ is good})\rightarrow 0.$$

The sets $A_\alpha$ depend on $\gamma$, but by selecting a
subsequence, we may assume that for each $\alpha\in \F^\ell$,
either \begin{equation} \label{small}
 |A_\alpha|=a_\alpha \mbox{ for some finite }
a_\alpha \end{equation} or
\begin{equation}\label{big}
|A_\alpha|\rightarrow \infty.
\end{equation}

Let $S=\{\alpha:\alpha \mbox{ is of type } \eqref{small} \}$. Let
$z$ be a random word as above (length  $n_\gamma$ with $r_\gamma
\leq (1/2-\varepsilon)n_\gamma$ non-zero coordinates). Let
$\mathcal{E}_1$ be the event that $\ma(z)=0$ for each $\alpha$
of type \eqref{small}. The bits $z(i)$ for the finitely many indices
$i\in A_\alpha$ for some $\alpha$ of type \eqref{small} are
asymptotically independent and each is 0 with probability
$(n_\gamma-r_\gamma)/n_\gamma>1/2$.

Hence
$$\liminf_{\gamma\rightarrow\infty} \P(\mathcal{E}_1)\geq
2^{-\sum_{\alpha\in S}a_\alpha}>0.$$
(This depends on $a_\alpha$, but we have chosen them and they are
now fixed). Given $\mathcal{E}_1$, for every  $\alpha\notin S$
(i.e., $\alpha$ is of type \eqref{big}) the random variable
$\ma(z)$ has a hypergeometric distribution with mean
$$\frac{r_\gamma}{n_\gamma-\sum_{\alpha\in S}a_\alpha }|A_\alpha|$$
and it follows by the law of large numbers that
$$\P\left(\left|\frac{\ma(z)}{|A_\alpha|}-\frac{r_\gamma}{n_\gamma}\right|
< \varepsilon \Bigm| \mathcal{E}_1\right)\rightarrow 1.$$
Since $r_\gamma/n_\gamma\leq 1/2-\varepsilon$, it follows that
$$\P\left(\frac{\ma(z)}{|A_\alpha|}
<\frac{1}{2}\Bigm| \mathcal{E}_1\right)\rightarrow 1$$ for each
$\alpha\notin S$.
 Hence, with probability
$(1+o(1))\P(\mathcal{E}_1)$,
$$
\begin{cases}
\ma(z)=0, & \alpha\in S,\\
\ma(z)<\frac{1}{2}|A_\alpha|, & \alpha\notin S,
\end{cases}
$$
and then $z$ is good.

Hence
$$\liminf_{\gamma\rightarrow \infty} \P(z \mbox{ is good})\geq
\liminf_{\gamma\rightarrow \infty}\P(\mathcal{E}_1)>0,$$ a
contradiction.
\end{proof}

For the remaining two proofs we will use the central limit theorem in
its simplest version, for symmetric binomial variables. (This was also
historically the first version, proved by de Moivre in 1733
\cite{deMoivre1733,Archibald}.)
We let, for $N\ge1$,
 $X_N$ denote a binomial random variable with the distribution
$\Bi(N,1/2)$. The central limit theorem says that
$(X_N-N/2)/\sqrt{N/4}$ converges in distribution to the standard
normal distribution $N(0,1)$, which means that if $Z\sim N(0,1)$, then
for any interval $I\subseteq\mathbb R$,
\begin{equation}
  \label{clt}
\P\Bigpar{\frac{X_N-N/2}{\sqrt{N/4}}\in I}
\to \P(Z\in I)
\qquad
\text{as }
N\to\infty.
\end{equation}
We will also need the more precise local central limit theorem which
says that if $x_N$ is any sequence of integers, then, as $N\to\infty$,
\begin{equation}
  \label{llt}
\P(X_N=x_N)
=\binom N{x_N}2^{-N}
=(2/\pi N)\qq\Bigpar{e^{-2(x_N-N/2)^2/N}+o(1)}.
\end{equation}
(This is a simple consequence of Stirling's formula.)

\begin{proof}[Proof of \refT{T2m}]
Let $\na=|A_\ga|$, and note that $\sum_{\ga\in\F^\ell} n_\ga=n$.
Fix an index $\gao$ with $n_\gao\ge n/2^{\ell}$ (for example the index
maximizing $\na$). Let $\ca=\set{\ga\in\F^\ell:\na>0}$ and
$\cax=\ca\setminus\set{\gao}$. Consider a random $z\in\F^n$. The
numbers $\ma=\ma(z)$ thus are independent binomial random variables:
$\ma\sim\Bi(\na,1/2)$.
Let $\gd=2^{-\ell-1}a$. Let $\ea$ be the event
\begin{equation}
  \label{ea}
\na/2 > \ma
\ge
\na/2-\ceil{\gd\sqrt n} ,
\end{equation}
for $\ga\in\cax$,  let $\eao$ be the event
\begin{equation}
  \label{eao}
\mao=r-\sumax\ma,
\end{equation}
and let $\ce=\EA$.
Assume in the sequel that $\sqrt n\ge 2^{\ell+1}/a$.
If $\ce$ holds, then
\begin{equation}
  \label{u3}
  \begin{split}
\mao=r-\sumax\ma
&\le n/2-a\sqrt n-\sumax(\na/2-\gd\sqrt{n} -1)
\\&
<n_\gao/2-a\sqrt n+2^\ell\gd\sqrt n+2^{\ell}
\le n_\gao/2,
  \end{split}
\end{equation}
and thus $z$ is good;
further, $d(z,x)=\sum_\ga\ma=r$.
It thus suffices to prove that $\P(\ce)\ge c n\qqw$,
since then
the number of good words
$z$ with $d(z,x)=r$
is at least $\P(\ce)2^n\ge c n\qqw 2^n$,
and further $\binom nr \le n\qqw 2^n$ by
\eqref{llt} (at least for large $n$).

First, let
\begin{equation*}
  p_N=\P(N/2 > X_N \ge N/2-\ceil{\gd\sqrt N}).
\end{equation*}
Note that
$p_N\ge \P(X_N =\floor{(N-1)/2})>0$ for every $N\ge1$, and that the
central limit theorem \eqref{clt} shows that as $N\to\infty$,
$p_N\to\P(0\ge Z\ge -2\gd)>0$.
Hence, $\px=\inf_{N\ge1} p_N>0$. Consequently, for $\ga\in\cax$,
$\P(\ea)\ge p_{\na}\ge\px$. Moreover, the events $\ea$, $\ga\in\cax$, are
independent, and thus
\begin{equation*}
  \P\Bigpar{\EAX}
=\prod_{\ga\in\cax} \P(\ea)
\ge \px^{2^\ell}.
\end{equation*}

Secondly, if \eqref{ea} holds for $\ga\in\cax$, then
$r-\sumax\ma<\nao/2$ by the calculation in \eqref{u3}, and
\begin{equation*}
  \begin{split}
r-\sumax\ma \ge n/2-b\sqrt n-\sumax \na/2 =\nao/2-b\sqrt n \ge
\nao/2-b2^{\ell/2}\sqrt {\nao}.
  \end{split}
\end{equation*}
 The random variable
$\mao$ is independent of \set{\ma:\ga\in\cax}, and
$\mao\sim\Bi(\nao,1/2)$. Thus, the local limit theorem \eqref{llt}
shows that for every set of numbers $\ma$, $\ga\in\cax$, satisfying
\eqref{ea},
\begin{equation*}
  \begin{split}
\P\bigpar{\eao\mid\ma,\,\ga\in\cax}
&=(2/\pi\nao)\qq\Bigpar{\exp\Bigpar{-2\Bigpar{r-\sumax\ma-\nao/2}^2/\nao}
+o(1)}
\\
&\ge
(2\nao)\qqw\Bigpar{\exp\Bigpar{-2^{\ell+1}b^2}+o(1)}
\ge \cc n\qqw
  \end{split}
\end{equation*}
for some $\ccx>0$, provided $n$, and thus also $\nao\ge2^{-\ell} n$, is
large enough.
Consequently, for large $n$,
\begin{equation*}
  \begin{split}
\P(\ce)
=\P\Bigpar{\EA}
=\P\Bigpar{\eao\mid\EAX}\P\Bigpar{\EAX}
\ge
\ccx n\qqw\px^{2^\ell}
=c n\qqw,
  \end{split}
\end{equation*}
which completes the proof.
\end{proof}

\begin{proof}[Proof of \refT{T3m}]
  Let $\na$, $\gao$, $\ca$ and $\cax$ be as in the preceding proof and
  consider again a random $z\in\F^n$.
Define the numbers $\ta$, $\ga\in\ca$, by
\begin{equation*}
  \ta=
  \begin{cases}
\floor{(\na-1)/2}, & \ga\in\cax,
\\
r-\sumax \ta, & \ga=\gao,
  \end{cases}
\end{equation*}
and let $\ce$ be the event
\begin{equation*}
  \ma=\ta,
\qquad \ga\in\ca.
\end{equation*}
Note that
\begin{equation*}
  \tao\le r-\sumax (\na/2-1)
\le
\nao/2-L+|\cax|<\nao/2,
\end{equation*}
and thus  $\ce$ implies that $z$ is good and
$d(z,x)=\sum_\ga\ma=r$.

Since also $\tao\ge r-\sumax\na/2\geq \nao/2-L-1$, it follows from
\eqref{llt} that for some
constant $\cc>0$ (depending on $L$) and every $n\ge 2^\ell L$,
\begin{equation*}
  \P(\ma=\ta)\ge\ccx \na\qqw\ge\ccx n\qqw
\end{equation*}
for every $\ga\in\ca$, and thus
\begin{equation*}
 \P(\ce)=\prod_{\ga\in\ca}\P(\ma=\ta)
\ge \cc n^{-2^\ell/2},
\end{equation*}
which completes the proof.
\end{proof}

\section{Construction of small $(r,\leq 2)$-identifying codes}
\label{Sconstr}

\begin{proof}[Proof of \refT{Tconstr}]
We make an explicit construction  for $r= \lfloor n/2\rfloor-1$. If
$2\leq n\leq 3$, then $r=0$ and we trivially may take $C=\F^n$.
Furthermore, the following few values are
known (see \cite{hlr,exoo}) $M_1^{(\leq 2)}(4)=11$, $M_1^{(\leq
2)}(5)=16$, $M_2^{(\leq 2)}(6)\leq 22$. So, we may assume that
$n\geq 7$.

Let  $C_0$ consist of the words $c_0=0\in \F^n$ and $c_i\in \F^n$ such
that $\supp(c_i)=\{1,2,\dots,i\}$ for $i=1,2,\dots,n$.
Clearly $|C_0|=n+1$.
Let $C_u=\{a\in \F^n\mid a\in C_0 \mbox{ or } \overline{a} \in
C_0\}$. Now $|C_u|=2n$.
The code which we claim to be $(r,\leq 2)$-identifying for
$r=\lfloor n/2 \rfloor -1$ is then the following
\begin{equation}\label{C}
  \begin{split}
C&=
\{ c\in \F^n \mid 1\leq d(c,a)\leq 2 \mbox{ for some } a\in C_u\}
\\
&=
\{ c\in \F^n \mid  d(c,a)= 2 \mbox{ for some } a\in C_u\},
  \end{split}
\end{equation}
where the equality follows since every word in $C_u$ has two
neighbours in $C_u$. Obviously, $|C|\leq \binom n2|C_u| =n^3-n^2$.
(We are interested in the order of growth, so this estimate is
enough for our purposes. However, with some effort one can check
that $|C|=n^3-5n^2+4n$ for $n\ge7$.)

We consider separately the cases $n$ even and $n$ odd.

1) Let first $n$ be odd.
The code $C_0$ is such that from every word $x\in
\F^n$ we have a codeword exactly at distance $(n-1)/2$. Indeed,
either $d(x,0)>(n-1)/2$ and $d(x,1)\leq (n-1)/2$ or $d(x,0)\leq
(n-1)/2$ and $d(x,1)>(n-1)/2$. Moving (in the first case --- the
second case is analogous) from the codeword $c_0=0$ to $c_n=1$
visiting every codeword $c_i$ $(i=1,\dots,n)$, there exists an index
$i$ such that $d(x,c_i)=(n-1)/2$, since every move between two
codewords $c_i$ and $c_{i+1}$ changes the distance by $\pm 1$.

Now we need to show that
$$I_r(X)\neq I_r(Y)$$
for any two distinct subsets $X\subseteq \F^n$ and $Y\subseteq \F^n$
where $|X|\leq 2$ and $|Y|\leq 2$. Assume to the contrary that
$I_r(X)= I_r(Y)$ for some $X,Y\subseteq \F^n$ with $|X|,|Y|\leq 2$
and $X\neq Y$.

Without loss of generality, we can assume that $|X|\geq |Y|$ and
that we have a word $x\in X\setminus Y$. Using the property of
$C_u$, we know that there exists a codeword $a\in C_u$ such that
$d(x,a)=(n-1)/2$ and $d(x,\overline{a})=(n+1)/2$. We concentrate on
the words in the sets $S_1(a)\cup S_2(a)$ and $S_1(\overline{a})\cup
S_2(\overline{a})$ which all belong to $C$. Since $I_r(x)\subseteq
I_r(X)$, we know that the sets
\begin{equation}\label{sets} I_r(X)\cap
S_1(a), \quad I_r(X)\cap S_2(a) \mbox { and } I_r(X)\cap
S_2(\overline{a})
\end{equation}
are all nonempty.
 By the symmetry of $\F^n$, we can
assume without loss of generality, that $a=0$ (and so,
$\overline{a}=1$).

Since $I_r(X)\cap S_1(0)$ is nonempty, there must be $\gamma\in Y$
such that $w(\gamma)\leq (n-1)/2$.

(i) Suppose first that $w(\gamma)\leq (n-5)/2$. This implies that
$S_1(0)\subseteq I_r(\gamma)\subseteq I_r(Y)=I_r(X)$.
Consequently, there must exist $y\in
X$ ($y\neq x$) such that $w(y)\leq (n-1)/2$, since $x$ does not
$r$-cover all of $S_1(0)$.
Since $|X|\le2$, thus $X=\{x,y\}$.

In order to cover the (nonempty) set $I_r(X)\cap S_2(1)$, there has
to be $\beta$ in $Y$ ($\beta \neq \gamma)$ such that $w(\beta)\geq
(n-1)/2$.
Thus $Y=\{\gamma,\beta\}$.
If $w(\beta)>(n-1)/2$, then $I_r(\beta)$ (and hence
$I_r(Y))$ contains elements from $S_1(1)$, but the set $I_r(X)\cap
S_1(1)$ is empty, immediately giving a contradiction. If
$w(\beta)=(n-1)/2$, then $I_r(x)$ contains a codeword not in
$I_r(Y)$. Indeed, since $x\neq \beta$ (and $w(x)=w(\beta)$), then
there exists an index $j\in \supp(\beta)$ such that $j\notin
\supp(x)$. This implies that the needed codeword, say $c'$, is found
in $S_2(1)$ by taking $\supp(\overline{c'})=\{i,j\}$ for any
$i\notin \supp(x)$ and $i\neq j$ --- clearly,
$\beta$ cannot $r$-cover this codeword
and $\gamma$ cannot $r$-cover any word in $S_2(1)$.

(ii) Assume then that $w(\gamma)=(n-3)/2$. Now $\gamma$ cannot
$r$-cover all the words in $I_r(x)\cap S_1(0)$,  so there must be
$\beta \in Y$ such that $w(\beta)\leq (n-1)/2$.
If $w(\beta)< (n-1)/2$, then $I_r(Y)\cap S_2(1)=\emptyset$
which contradicts $I_r(X)\cap S_2(1)\neq\emptyset$.
If $w(\beta)=(n-1)/2$,
we are done as in (i),
using again $x\neq \beta$.

(iii) Let then $w(\gamma)=(n-1)/2$. Because there are codewords in
$I_r(x)\cap S_1(0)$ which are not $r$-covered by $\gamma$, it
follows that
there exists $\beta\in Y$ with
$w(\beta)\le(n-1)/2$.
By the previous cases, it suffices
to consider
$w(\beta)=(n-1)/2$, since otherwise we can interchange $\beta$ and $\gamma$.
Let $i\in \supp(x)$ be such that $i\notin \supp(\gamma)$ and
$j\in\supp(x)$ such that $j\notin \supp(\beta)$. Since $x\notin Y$,
such indices (it is possible that $i=j$) exist. When $i\neq j$, a
codeword $c\in S_2(0)$ such that  $\supp( c)=\{i,j\}$, gives a
contradiction. If $i=j$, then we
pick a codeword with $\supp(c)=\{i,k\}$ where $k\in \supp(x)$,
$i\neq k$.

\medskip

2) Let now $n$ be even. Take $C_u$ as in the odd case; it has now
the analogous property that from every word $x\in \F^n$ there is a
codeword $a\in C_u$ such that $d(x,a)=d(x,\overline{a})=n/2$. Let
$C$ be defined also as above. We will show that it is $(r,\leq
2)$-identifying for $r=n/2-1$. If $I_r(X)=I_r(Y)$, we can again
assume that $|X|\geq |Y|$ and choose  $x\in X\setminus Y$. We know
that there is $a\in C_u$ such that $d(x,a)=n/2$. The sets
\eqref{sets} as well as now the set $I_r(x)\cap S_1(\overline{a})$
are nonempty. Again it suffices to consider $a=0$. Since $I_r(x)\cap
S_1(0)$ is nonempty, so there must be a word $\gamma\in Y$ such that
$w(\gamma)\leq n/2$.

(i) Suppose first that $w(\gamma)\leq n/2-2$. Then $S_1(0)\subseteq
I_r(\gamma)$. Since $S_1(0)\nsubseteq I_r(x)$, this implies that
there is $y\in X$, $y\neq x$, such that $w(y)\leq n/2$.

Let first $w(y)\leq n/2-1$. Subsequently, neither $y$ nor $\gamma$
$r$-covers any of the codewords of $S_1(1)$ whereas $|I_r(x)\cap
S_1(1)|=n/2$. Hence there has to be $\beta\in Y$ such that
$I_r(x)\cap S_1(1)=I_r(\beta)\cap S_1(1)$. However, this implies
that $x=\beta$, a contradiction.

Assume next that $w(y)=n/2$. Due to the fact that $\gamma$
$r$-covers all the codewords $S_1(0)$ we know that $y=\overline{x}$.
Hence $S_1(1)\subseteq I_r(X)=I_r(Y)$.
On the other hand,
$S_1(1)\cap I_r(\gamma)=\emptyset$.
Consequently,
$S_1(1)\subseteq I_r(\beta)$ which implies
$w(\beta)\geq n/2+2$.
Thus $\beta$ does not $r$-cover any of
the words in $S_2(0)$.
\begin{itemize}
\item If $w(\gamma)\leq n/2-3$, then $\gamma$ $r$-covers all of
$S_2(0)$. However, all of $S_2(0)$ is not contained in $I_r(X)$.
Indeed, take $i\in \supp(x)$ and $j\notin \supp(x)$ (notice that now
$X=\{x,\overline{x}\})$. The codeword $c'$ of $S_2(0)$ which has
$\supp(c')=\{i,j\}$ does not belong to $I_r(X)$.

\item If $w(\gamma)= n/2-2$, then $\supp(x)$  has (at least) two distinct indices, say $i$ and $j$,
 which are not in $\supp(\gamma)$. Consequently, the codeword $c'$ in $S_2(0)$ with $\supp(c')=\{i,j\}$ belongs to $I_r(x)$ but not to
 $I_r(Y)$, a contradiction.

\end{itemize}

(ii) Let now $w(\gamma)=n/2-1$. Now $I_r(x)\cap S_1(0)$ has one more
codeword than $I_r(\gamma)\cap S_1(0)$ and, therefore,
there exists $\beta\in Y$ such that
$w(\beta)\leq
n/2$ (or we are done). But now $I_r(x)\cap S_1(1)$ contains at least
one codeword not in $I_r(\beta)$ --- notice that $\gamma$ does not
$r$-cover any words in $S_1(1)$.

(iii) Assume finally that $w(\gamma)=n/2$.
Since $\gamma\neq x$, there must exist an index
$i\in \supp(x)$ such that $i\notin \supp(\gamma)$. Consequently, the
$r$ codewords $z_j\in  S_2(0)\cap I_r(x)$ such that
$\supp(z_j)=\{i,j\}$, where $j\in\supp(x)$ and $j\neq i$, do not
belong to $I_r(\gamma)$ and hence must belong to $I_r(\beta)$
for some other $\beta\in Y$
or we
are done.
By the previous cases and
symmetry with respect to $a$ and $\overline{a}$, we can also assume
that $w(\beta)=n/2$.
However, this means that $\supp(z_j)\subseteq
\supp(\beta)$ for all $j$. Subsequently, $x=\beta$ or $w(\beta)>
n/2$ and we get a contradiction in both cases.
\end{proof}


\begin{ack}
  This research was done mainly during the
Workshop on Codes and Discrete Probability in Grenoble, 2007.
\end{ack}

\end{document}